\overfullrule=0pt
\centerline {\bf Another multiplicity result for the periodic solutions of certain systems}\par
\bigskip
\bigskip
\centerline {BIAGIO RICCERI}\par
\bigskip
\bigskip
{\bf Abstract:} In this paper, we deal with a problem of the type
$$\cases{(\phi(u'))'=\nabla_xF(t,u) & in $[0,T]$\cr & \cr
u(0)=u(T)\ , \hskip 3pt u'(0)=u'(T)\ ,\cr}$$
where, in particular, $\phi$ is a homeomorphism from an open ball of ${\bf R}^n$ onto ${\bf R}^n$. Using the theory developed by
Brezis and Mawhin in [1] jointly with our minimax theorem proved in [3], we obtain a general multiplicity result, under assumptions of qualitative
nature only. Three remarkable corollaries are also presented.\par
\bigskip
\bigskip
\bigskip
\bigskip
{\bf 1. Introduction}\par
\bigskip
In what follows, $L, T$ are two fixed positive numbers. For each $r>0$, we set $B_r=\{x\in {\bf R}^n :|x|<r\}$ ($|\cdot|$ being the Euclidean norm on ${\bf R}^n$)
and $\overline {B_r}$ is the closure of $B_r$.\par
\smallskip
We denote by ${\cal A}$ the family of all homeomorphisms $\phi$ from $B_L$ onto ${\bf R}^n$ such that $\phi(0)=0$ and $\phi=\nabla\Phi$,
where the function $\Phi:\overline {B_L}\to ]-\infty,0]$ is continuous and strictly convex in $\overline {B_L}$, and of class $C^1$ in $B_L$.
Notice that $0$ is the unique global minimum of $\Phi$ in $\overline {B_L}$.
\par
\smallskip
We denote by ${\cal B}$ the family of all functions $F:[0,T]\times {\bf R}^n\to {\bf R}$ which are measurable in $[0,T]$, of class $C^1$ in
${\bf R}^n$ and such that $\nabla_xF$ is measurable in $[0,T]$ and, 
for each $r>0$, one has $\sup_{x\in B_r}|\nabla_x F(\cdot,x)|\in L^1([0,T])$, with $F(\cdot,0)\in L^1([0,T])$.\par
\smallskip
Given $\phi\in {\cal A}$ and $F\in {\cal B}$, we consider the problem
$$\cases{(\phi(u'))'=\nabla_xF(t,u) & in $[0,T]$\cr & \cr
u(0)=u(T)\ , \hskip 3pt u'(0)=u'(T)\ .\cr}\eqno{(P_{\phi,F})}$$
A solution of this problem is any function $u:[0,T]\to {\bf R}^n$ of class $C^1$, with $u'([0,T])\subset
B_L$, $u(0)=u(T)$, $u'(0)=u'(T)$,  such that the composite function $\phi\circ u'$ is absolutely
continuous in $[0,T]$ and one has $(\phi\circ u')'(t)=\nabla_xF(t,u(t))$ for a.e. $t\in [0,T]$.\par
\smallskip
Now, we set
$$K=\{u\in \hbox {\rm Lip}([0,T],{\bf R}^n) : |u'(t)|\leq L\hskip 5pt for\hskip 3pt a.e.\hskip 3pt t\in [0,T] , u(0)=u(T)\}\ ,$$
Lip$([0,T],{\bf R}^n)$ being the space of all Lipschitzian functions from $[0,T]$ into ${\bf R}^n$.\par
\smallskip
Clearly, one has
$$\sup_{[0,T]}|u|\leq LT + \inf_{[0,T]}|u| \eqno{(1.1)}$$
for all $u\in K$. 
\smallskip
Next, consider the functional $I:K\to {\bf R}$ defined by
$$I(u)=\int_0^T(\Phi(u'(t))+F(t,u(t)))dt$$
for all $u\in K$.\par
\smallskip
In [1], Brezis and Mawhin proved the following result:\par
\medskip
THEOREM 1.A ([1], Theorem 5.2). - {\it Any global minimum of $I$ in $K$ is a solution of problem $(P_{\phi,F})$.}\par
\medskip
In [4], using Theorem 1.A jointly with the theory developed in [2], we obtained the following multiplicity theorem:\par
\medskip
THEOREM 1.B ([4], Theorem 3.1). - {\it Let $\phi\in {\cal A}$, $F\in {\cal B}$ and $G\in C^{1}({\bf R}^n)$.
  Moreover, let $\gamma:[0,+\infty[\to {\bf R}$ be a 
convex strictly increasing function such that $\lim_{s\to +\infty}{{\gamma(s)}\over {s}}=+\infty$.
Assume that the following assumptions are
satisfied:\par
\noindent
$(i_1)$\hskip 3pt 
for a.e. $t\in [0,T]$ and for every $x\in {\bf R}^n$, one has
$$\gamma(|x|)\leq F(t,x)\ ;$$
\noindent
$(i_2)$\hskip 3pt $\liminf_{|x|\to +\infty}{{G(x)}\over {|x|}}>-\infty\ ;$\par
\noindent
$(i_3)$\hskip 3pt the function $G$ has no global minima in ${\bf R}^n$\ ;\par
\noindent
$(i_4)$\hskip 3pt  there exist two distinct points $x_1, x_2\in {\bf R}^n$ such that
$$\inf_{x\in {\bf R}^n}\int_0^TF(t,x)dt<\max\left \{ \int_0^TF(t,x_1)dt, \int_0^TF(t,x_2)dt\right\}$$
and
$$G(x_1)=G(x_2)=\inf_{B_c}G$$
where
$$c=LT+\gamma^{-1}\left ( {{1}\over {T}}\max\left \{ \int_0^TF(t,x_1)dt, \int_0^TF(t,x_2)dt\right\}\right )\ .$$
Then, for every $\psi\in L^{1}([0,T])\setminus \{0\}$, with $\psi\geq 0$,
there exists $\tilde\lambda>0$ such that the problem
$$\cases{(\phi(u'))'=\nabla_x(F(t,u)+\tilde\lambda \psi(t)G(u)) & in $[0,T]$\cr & \cr
u(0)=u(T)\ ,\hskip 3pt u'(0)=u'(T)\cr}$$
has at least two solutions which are global minima in $K$ of the functional 
$$u\to \int_0^T(\Phi(u'(t))+F(t,u(t))+\tilde\lambda \psi(t)G(u(t)))dt\ .$$}\par
\medskip
Clearly, condition $(i_4)$ is a little involved and has a typical quantitative nature, due to the presence of the constant
$c$. This kind of drawback, however, is largely compensated by the great generality of the conclusion, due to its validity for any
 $\psi\in L^{1}([0,T])\setminus \{0\}$, with $\psi\geq 0$.\par
\smallskip
The aim of the present short paper is to give a further contribution to the subject, adopting assumptions of qualitative nature only, in the spirit
of Theorem 1.2 of [5].\par
\bigskip
{\bf 2. Results}\par
\bigskip
The two main tools we will use to prove our main theorem are as follows:\par
\medskip
THEOREM 2.A ([3], Theorem 1.2). - {\it Let $X$ be a topological space, $E$ a real vector
space, $Y\subseteq E$ a non-empty convex set and $J:X\times Y\to {\bf R}$ a function which is lower semicontinuous and inf-compact in $X$, 
and concave in $Y$. Moreover, assume that
$$\sup_Y\inf_XJ<\inf_X\sup_YJ\ .$$
Then, there exists $\hat y\in Y$ such that the function $J(\cdot,\hat y)$ has at least two global
minima.}
\medskip
PROPOSITION 2.A ([5], Proposition 2.2). - {\it Let $X, Y$ be two non-empty sets and $f:X\to {\bf R}$, $g:X\times Y\to {\bf R}$ two given
functions. Assume that there are two sets $A, B\subset X$ such that:\par
\noindent
$(a)$\hskip 5pt $\sup_Af<\inf_Bf$\ ;\par
\noindent
$(b)$\hskip 5pt $\sup_{y\in Y}\inf_{x\in A}g(x,y)\leq 0$\ ;\par
\noindent
$(c)$\hskip 5pt $\inf_{x\in B}\sup_{y\in Y}g(x,y)\geq 0$\ ;\par
\noindent
$(d)$\hskip 5pt $\inf_{x\in X\setminus B}\sup_{y\in Y}g(x,y)=+\infty$\ .\par
Then, one has
$$\sup_{y\in Y}\inf_{x\in X}(f(x)+g(x,y))<\inf_{x\in X}\sup_{y\in Y}(f(x)+g(x,y))\ .$$}\par
\medskip
A set $Y\subseteq L^1([0,T])$ is said to have property $(P)$ if 
$$\sup_{\psi\in Y}\int_0^T\psi(t)h(t)dt=+\infty$$
 for all $h\in C^0([0,T])\setminus \{0\}$.\par
\smallskip
Our main result is as follows:\par
\medskip
THEOREM 2.1. - {\it Let $\phi\in {\cal A}$, $F\in {\cal B}$ and $G\in C^1({\bf R}^n)$. Assume that:\par
\noindent
$(a_1)$\hskip 5pt there exists $q>0$ such that
$$\lim_{|x|\to +\infty}{{\inf_{t\in [0,T]}F(t,x)}\over {|x|^q}}=+\infty$$
and
$$\limsup_{|x|\to +\infty}{{|G(x)|}\over {|x|^q}}<+\infty\ ;$$
\noindent
$(a_2)$\hskip 5pt there exists $r\in \left ] \inf_{{\bf R}^n}G,\sup_{{\bf R}^n}G\right [$ such that
$$\max \left \{ \inf_{x\in G^{-1}(]-\infty,r[)}\int_0^TF(t,x)dt, \inf_{x\in G^{-1}(]r,+\infty[)}\int_0^TF(t,x)dt\right \}<
\int_0^T\inf_{x\in G^{-1}(r)} F(t,x)dt\ .$$
Then, for every non-empty convex set $Y\subseteq L^{\infty}([0,T])$ with property $(P)$, there exists $\psi\in Y$ such that the problem
$$\cases{(\phi(u'))'=\nabla_x(F(t,u)+\psi(t)G(u)) & in $[0,T]$\cr & \cr
u(0)=u(T)\ ,\hskip 3pt u'(0)=u'(T)\cr}$$
has at least two solutions which are global minima in $K$ of the functional 
$$u\to \int_0^T(\Phi(u'(t))+F(t,u(t))+\psi(t)G(u(t)))dt\ .$$}\par
\smallskip
PROOF. Fix a non-empty convex set $Y\subseteq L^{\infty}([0,T])$ with property $(P)$.
Let $C^0([0,T], {\bf R}^n)$ be the space of all continuous functions from $[0,T]$ into ${\bf R}^n$, with the norm $\sup_{[0,T]}|u|$.
To achieve the conclusion, we are going to apply Theorem 2.A taking $X=K$, regarded as a subset of $C^0([0,T], {\bf R}^n)$ with the relative topology,
and $J:K\times Y\to {\bf R}$ defined by
$$J(u,\psi)=\int_0^T(\Phi(u'(t))+F(t,u(t)))dt+\int_0^T\psi(t)(G(u(t))-r)dt$$
for all $(u,\psi)\in K\times Y$. Clearly, $J(u,\cdot)$ is concave in $Y$. Fix $\psi\in Y$. 
By Lemma 4.1 of [1], $J(\cdot,\psi)$ is lower semicontinuous in $K$. Let us show that $J(\cdot,\psi)$ is inf-compact in $K$.
By $(a_1)$, there exist $k, \delta, \nu>0$, with 
$$\nu>k\|\psi\|_{L^{\infty}}\ , \eqno{(2.1)}$$
such that
$$|G(x)|\leq k(|x|^q+1) \eqno{(2.2)}$$
for all $x\in {\bf R}^n$ and
$$F(t,x)\geq \nu|x|^q$$
for all $t\in [0,T]$ and $x\in {\bf R}^n\setminus B_{\delta}$. Since $F\in {\cal B}$, there exists $M\in L^1([0,T])$ such that
$$|\nabla_x F(t,x)|\leq M(t)$$
for all $t\in [0,T]$ and $x\in B_{\delta}$. By the mean value theorem, we have
$$F(t,x)-F(t,0)=\langle \nabla_x F(t,\xi),x\rangle$$
for some $\xi$ in the segment joining $0$ and $x$. Consequently, for all $t\in [0,T]$ and $x\in B_{\delta}$, we have
$$||F(t,x)|-|F(t,0)||\leq |F(t,x)-F(t,0)|\leq \delta M(t)$$
and so, if we put
$$\beta(t)=\nu\delta^q+M(t)\delta+|F(t,0)|\ ,$$
we have
$$F(t,x)\geq \nu|x|^q-\beta(t) \eqno{(2.3)}$$
for all $t\in [0,T]$ and $x\in {\bf R}^n$. Now, set
$$\eta=-\int_0^T\beta(t)dt+\Phi(0)T-r\int_0^T|\psi(t)|dt\ ,$$
and
$$\eta_1=\eta-k\int_0^T|\psi(t)|dt\ .$$
For each $u\in K$, with $\sup_{[0,T]}|u|\geq LT$,  taking $(1.1)$, $(2.1)$, $(2.2)$, $(2.3)$ into account, we  have
$$J(u,\psi)\geq \int_0^TF(t,u(t))dt-\int_0^T|\psi(t)G(u(t))|dt+\Phi(0)T-r\int_0^T\psi(t)dt$$
$$\geq\nu\int_0^T|u(t)|^qdt-\int_0^T|\psi(t)G(u(t))|dt+\eta\geq \nu\int_0^T|u(t)|^qdt-k\|\psi\|_{L^{\infty}}\int_0^T|u(t)|^qdt+\eta_1$$
$$\geq (\nu- k\|\psi\|_{L^{\infty}})T\left (\inf_{[0,T]}|u|\right )^q+\eta_1\geq 
(\nu- k\|\psi\|_{L^{\infty}})T\left (\sup_{[0,T]}|u|-LT\right )^q+\eta_1$$
Consequently
$$\sup_{[0,T]}|u|\leq \left ( {{J(u,\psi)-\eta_1}\over {(\nu-k\|\psi\|_{L^{\infty}})T}}\right )^{1\over q}+LT\ .\eqno{(2.4)}$$
Fix $\rho\in {\bf R}$. By $(2.4)$, the set $C_{\rho}:=\{u\in K : J(u,\psi)\leq \rho\}$ turns out to be bounded. 
Moreover, the functions belonging to $C_{\rho}$ are
equi-continuous since they lie in $K$. As a consequence, by the Ascoli-Arzel\`a theorem, $C_{\rho}$ is relatively compact in $C^0([0,T],{\bf R}^n)$.  By
lower semicontinuity, $C_{\rho}$ is closed in $K$. But $K$ is closed in $C^0([0,T], {\bf R}^n)$ and hence $C_{\rho}$ is compact. The inf-compactness of
$J(\cdot,\psi)$ is so shown. 
Now, to obtain the strict minimax inequality required by Theorem 2.1, we use Proposition 2.A. By $(a_2)$, there are $x_1, x_2\in {\bf R}^n$ such that
$$G(x_1)<r<G(x_2) \eqno{(2.5)}$$
and
$$\max\left \{ \int_0^TF(t,x_1)dt, \int_0^TF(t,x_2)dt\right \}<\int_0^T\inf_{x\in G^{-1}(r)} F(t,x)\ .\eqno{(2.6)}$$
Now, put
$$A=\{x_1, x_2\}\ ,$$
$$B=\{u\in K : u([0,T])\subseteq G^{-1}(r)\}$$
and define $f:K\to {\bf R}$, $g:K\times Y\to {\bf R}$ by
$$f(u)=\int_0^T(\Phi(u'(t))+F(t,u(t)))dt\ ,$$
$$g(u,\psi)=\int_0^T\psi(t)(G(u(t))-r)dt$$
for all $u\in K$, $\psi\in Y$. Since the constant functions (from $[0,T]$ into ${\bf R}^n$) belong to $K$, we think of $A$ as a subset of $K$.
With these choices, in connection with Proposition 2.A, $(a)$ is a simple consequence of $(2.6)$; $(b)$ follows immediately from $(2.5)$; $(c)$ is obvious since
$g(u,\psi)=0$ for all $u\in B$. Finally, concerning $(d)$, observe that, if $u\in K\setminus B$, then the continuous function $G\circ u-r$ is not zero
and hence $\sup_{\psi\in Y}g(u,\psi)=+\infty$ since $Y$ has property $(P)$. Therefore, Propostion 2.A ensures that
$$\sup_Y\inf_KJ<\inf_K\sup_YJ\ .$$
Now, our conclusion follows directly from Theorem 2.1 and Theorem 1.A.\hfill $\bigtriangleup$\par
\medskip
We now point out three remarkable corollaries of Theorem 2.1.\par
\medskip
COROLLARY 2.1. - {\it Let $\phi\in {\cal A}$, $F\in {\cal B}$ and $G\in C^1({\bf R}^n)$.
Besides condition $(a_1)$, assume that
$F(t,\cdot)$ is even for all $t\in [0,T]$, that $G$ is odd and that
$$\inf_{x\in {\bf R}^n}\int_0^TF(t,x)dt<\int_0^T\inf_{x\in G^{-1}(0)} F(t,x)dt\ .$$
Then, the conclusion of Theorem 2.1 holds.}\par
\smallskip
PROOF. By assumption, there is $\hat x\in {\bf R}^n$ such that 
$$\int_0^TF(t,\hat x)dt<\int_0^T\inf_{x\in G^{-1}(0)} F(t,x)dt\ .$$
So, $G(\hat x)\neq 0$. Assume, for instance, that $G(\hat x)>0$. Then, since $G$ is odd, $G(-\hat x)<0$. But, since $F(t,\cdot)$ is even, we have
$$\int_0^TF(t,\hat x)dt=\int_0^TF(t,-\hat x)dt\ .$$
Therefore, condition $(a_2)$ is satisfied with $r=0$, and the conclusion of Theorem 2.1 follows.\hfill $\bigtriangleup$\par
\medskip
COROLLARY 2.2. - {\it Let  $\phi\in {\cal A}$ and let $F, G\in C^1({\bf R}^n)$, with $\lim_{|x|\to +\infty}G(x)=+\infty$.
Besides condition $(a_1)$, assume that
there exists a point $x_0\in {\bf R}^n$ which is, at the same time, the unique global minimum of $G$ and a strict local, not global,
minimum of $F$.\par
Then, for each $\gamma\in L^1([0,T])$, with $\inf_{[0,T]}\gamma>0$, and for each non-empty convex set $Y\subseteq L^{\infty}[0,T])$ with property $(P)$,
there exists $\psi\in Y$ such that the problem
$$\cases{(\phi(u'))'=\nabla_x(\gamma(t)F(u)+\psi(t)G(u)) & in $[0,T]$\cr & \cr
u(0)=u(T)\ ,\hskip 3pt u'(0)=u'(T)\cr}$$
has at least two solutions which are global minima in $K$ of the functional 
$$u\to \int_0^T(\Phi(u'(t))+\gamma(t)F(u(t))+\psi(t)G(u(t)))dt\ .$$}\par
\smallskip
PROOF. By assumption, there are $x_1\in {\bf R}^n$ and $\rho>0$ such that 
$$F(x_1)<F(x_0)<F(x) \eqno{(2.7)}$$ 
for all $x\in B(x_0,\rho)\setminus \{x_0\}$. Now, observe that,
since $G$ is inf-compact (being coercive) and $x_0$ is the unique global minimum of $G$, for each sequence $\{y_k\}$ in ${\bf R}^n$
such that $\lim_{k\to \infty}G(y_k)=G(x_0)$, we have $\lim_{k\to \infty}y_k=x_0$. As a consequence, we can fix $r>G(x_0)$ so that
$$G^{-1}(]-\infty,r])\subseteq B(x_0,\rho)\ .\eqno{(2.8)}$$
From $(2.7)$ and $(2.8)$, it follows that 
$$G(x_1)>r$$ as well as, by compactness,
$$F(x_0)<\inf_{x\in G^{-1}(r)}F(x)\ .$$
At this point, it is clear that, for each $\gamma\in L^1([0,T]$, with $\inf_{[0,T]}\gamma>0$, the function $(t,x)\to \gamma(t)F(x)$ satisfies
conditions $(a_1)$ and $(a_2)$ and the conclusion follows.\hfill $\bigtriangleup$
\medskip
Recall that a real-valued function on a convex subset of a vector space is said to be quasi-convex if its sub-level sets are convex.\par
\medskip 
COROLLARY 2.3. - {\it Let $n=1$ and let $\phi\in {\cal A}$, $F\in {\cal B}$, $G\in C^1({\bf R})$. Besides condition $(a_1)$,
assume that $G$ is strictly monotone and that $x\to \int_0^TF(t,x)dt$ is not quasi-convex.\par
Then, the conclusion of Theorem 2.1 holds.}\par
\smallskip
PROOF.  By assumption,
 there are $x_1, x_2, x_3\in {\bf R}$, with $x_1<x_3<x_2$, such that
$$\max\left \{ \int_0^TF(t,x_1)dt, \int_0^TF(t,x_2)dt\right \}<\int_0^TF(t,x_3)dt\ .$$
Moreover, the numbers $G(x_1)-G(x_3)$ and $G(x_2)-G(x_3)$ have opposite signs and $G^{-1}(G(x_3))=\{x_3\}$.
Therefore, condition $(a_2)$ is satisfied with $r=G(x_3)$, and the conclusion follows.\hfill $\bigtriangleup$
\bigskip
\bigskip
{\bf Acknowledgement.} The author has been supported by the Gruppo Nazionale per l'Analisi Matematica, la Probabilit\`a e 
le loro Applicazioni (GNAMPA) of the Istituto Nazionale di Alta Matematica (INdAM) and by the Universit\`a degli Studi di Catania, ``Piano della Ricerca 2016/2018 Linea di intervento 2". 
\vfill\eject
\centerline {\bf References}\par
\bigskip
\bigskip
\noindent
[1]\hskip 5pt H. BREZIS and J. MAWHIN, {\it Periodic solutions of Lagrangian systems of relativistic oscillators},
Commun. Appl. Anal., {\bf 15} (2011), 235-250.\par
\medskip
\noindent
[2]\hskip 5pt B. RICCERI, {\it Well-posedness of constrained minimization problems via saddle-points}, J. Global Optim., {\bf 40} (2008),
389-397.\par
\medskip
\noindent
[3]\hskip 5pt B. RICCERI, {\it On a minimax theorem: an improvement, a new proof and an overview of its applications},
Minimax Theory Appl., {\bf 2} (2017), 99-152.\par
\medskip
\noindent
[4]\hskip 5pt B. RICCERI, {\it Multiple periodic solutions of Lagrangian systems of relativistic oscillators}, in  ``Current Research in Nonlinear Analysis - In
Honor of Haim Brezis and Louis Nirenberg", Th. M. Rassias ed., 249-258, Springer, 2018.\par
\medskip
\noindent
[5]\hskip 5pt B. RICCERI, {\it Miscellaneous applications of certain minimax theorems II}, Acta Math. Vietnam., to appear.\par
\bigskip
\bigskip
\bigskip
\bigskip
Department of Mathematics and Informatics\par
University of Catania\par
Viale A. Doria 6\par
95125 Catania, Italy\par
{\it e-mail address}: ricceri@dmi.unict.it

\bye